\documentclass{article}
\usepackage{amssymb}
\usepackage{amsfonts}

\newtheorem{theorem}{Theorem}

\newtheorem{conjecture}[theorem]{Conjecture}
\newtheorem{corollary}[theorem]{Corollary}

\newtheorem{lemma}[theorem]{Lemma}

\newenvironment{proof}[1][Proof]{\noindent\textbf{#1.} }{\ \rule{0.5em}{0.5em}}
\input{tcilatex}
\begin{document}

\author{Zafer \c{S}iar and Refik Keskin \\
Bing\"{o}l University, Department of Mathematics, Bing\"{o}l/TURKEY\\
zsiar@bingol.edu.tr\\
Sakarya University, Department of Mathematics, Sakarya/TURKEY\\
rkeskin@sakarya.edu.tr}
\title{On the Exponential Diophantine Equation $(a^{n}-2)(b^{n}-2)=x^{2}$}
\maketitle

\begin{abstract}
In this paper, we consider the equation%
\begin{equation}
(a^{n}-2^{m})(b^{n}-2^{m})=x^{2},~x,n,m\in 
\mathbb{N}
.  \label{A}
\end{equation}%
By assuming the $abc$ conjecture is true, in \cite{Luca1}, Luca and Walsh
gave a theorem, which implies that the equation (\ref{A}) has only finitely
many solutions $n,x$ if $a$ and $b$ are different fixed positive integers.
We solve (\ref{A}) when $m=1$ and $(a,b)\in \left\{
(2,10),(4,100),(10,58),(3,45)\right\} .$ Moreover, we show that $%
(a^{n}-2)(b^{n}-2)=x^{2}$ has no solution $n,x$ if $2|n$ and $\gcd (a,b)=1.$
We also give a conjecture which says that the equation $%
(2^{n}-2)((2P_{k})^{n}-2)=x^{2}$ has only the solution $(n,x)=(2,Q_{k}),$
where $k>3$ is odd and $P_{k},Q_{k}$ are Pell and Pell Lucas numbers,
respectively. We also conjecture that if the equation $%
(a^{n}-2)(b^{n}-2)=x^{2}$ has a solution $n,x,$ then $n\leq 6$, where $%
2<a<b. $
\end{abstract}



\section{\protect\bigskip Introduction}

Over the last decades, several authors have dealt with \ the equation 
\begin{equation}
(a^{n}-1)(b^{n}-1)=x^{2},~x,n\in 
\mathbb{N}
,  \label{1}
\end{equation}%
where $a>1$ and $b>1$ are different fixed integers. Firstly, Szalay \cite%
{szly} handled this equation for $(a,b)=(2,3),(2,5),$ and $(2,2^{k}).$ He
showed that there is no solution if $(a,b)=(2,3),$ only the solution $n=1$
for $(a,b)=(2,5),$ and only the solution $k=2$ and $n=3$ for $(a,b)=$ $%
(2,2^{k})$ with $k>1.$ Then, in \cite{szly2}, the authors determined that
the equation (\ref{1}) has no solutions for $(a,b)=(2,6)$ and that has the
only solutions $(a,n,k)=(2,3,2),(3,1,5),$ and $(7,1,4)$ for $(a,b)=(a,a^{k})$
with $kn>2.$ Last result was extended by Cohn \cite{chn} to the case $%
a^{l}=b^{k}.$ He also proved that the equation (\ref{1}) has no solutions if 
$4|n$ except for $(a,b)=(13,239)$, in which case $n=4.$ Later, in \cite%
{LS,ZL,MT,GX}, the authors studied the equation (\ref{1}) for the different
values of $a$ and $b.$ Lastly, in \cite{rfk2}, the Keskin proved that the
equation (\ref{1}) has no solutions for $n>4$ with $2|n$ if $a$ and $b$ have
opposite parity. Keskin also proved that if $\gcd (a,b)=1,$ $2||n,$ and $%
n>4, $ then the equation (\ref{1}) has no solutions.

Motivated by the above studies, in this paper, we consider the equation%
\begin{equation}
(a^{n}-2^{m})(b^{n}-2^{m})=x^{2},~x,n,m\in 
\mathbb{N}
.  \label{3.1}
\end{equation}%
By assuming the $abc$ conjecture is true, in \cite{Luca1}, Luca and Walsh
gave a theorem, which implies that the equation (\ref{3.1}) has only
finitely many solutions $n,x$ if $a$ and $b$ are different fixed positive
integers. For more on the abc conjecture, one can consult \cite{Pa}. We
solve (\ref{3.1}) when $m=1$ and $(a,b)\in \left\{
(2,10),(4,100),(10,58),(3,45)\right\} .$ Moreover, we show that $%
(a^{n}-2)(b^{n}-2)=x^{2}$ has no solution $n,x$ if $2|n$ and $\gcd (a,b)=1.$
We also give a conjecture which says that the equation $%
(2^{n}-2)((2P_{k})^{n}-2)=x^{2}$ has only the solution $(n,x)=(2,Q_{k}),$
where $k>3$ is odd and $P_{k},Q_{k}$ are Pell and Pell Lucas numbers,
respectively. We also conjecture that if the equation $%
(a^{n}-2)(b^{n}-2)=x^{2}$ has a solution $n,x,$ then $n\leq 6$, where $%
2<a<b. $

\section{Preliminaries}

In this study, while solving the equation (\ref{3.1}), the first and second
kind of Lucas sequence $U_{n}(P,Q)$ and $V_{n}(P,Q)$ play an essential role.
So, we need to recall them.

Let $P$ and $Q$ be nonzero relatively prime integers such that $P^{2}+4Q>0.$
Define 
\[
U_{0}(P,Q)=0,~U_{1}(P,Q)=1,~U_{n+1}(P,Q)=PU_{n}(P,Q)+QU_{n-1}(P,Q)\text{ for 
}n\geq 1, 
\]%
\[
V_{0}(P,Q)=2,~V_{1}(P,Q)=P,~V_{n+1}(P,Q)=PV_{n}(P,Q)+QV_{n-1}(P,Q)\text{ for 
}n\geq 1. 
\]%
Sometimes, we write $U_{n}$ and $V_{n}$ instead of $U_{n}(P,Q)$ and $%
V_{n}(P,Q),$ respectively. For $(P,Q)=$ $(2,1),$ we have Pell and Pell-Lucas
sequences, $(P_{n})$ and $(Q_{n}).$ The following identities concerning the
sequences $\left( U_{n}\right) $ and $\left( V_{n}\right) ,$ which will be
used in the next section, are well known (see \cite{rbn,zfr2,Zfr}).

Let $d=(m,n).$ Then 
\begin{equation}
(U_{m},U_{n})=U_{d}  \label{0.4}
\end{equation}%
and 
\begin{equation}
(V_{m},V_{n})=\left\{ 
\begin{array}{c}
V_{d}\text{ \ if }m/d\text{ and }n/d\text{ are odd,} \\ 
1\text{ or }2\text{ \ \ \ otherwise.}%
\end{array}%
\right.  \label{0.5}
\end{equation}%
\begin{equation}
V_{2n}=V_{n}^{2}-2(-Q)^{n}  \label{0.11}
\end{equation}%
\begin{equation}
V_{3n}=V_{n}(V_{n}^{2}-3(-Q)^{n})  \label{0.9}
\end{equation}%
\begin{equation}
V_{n}(P,-1)=U_{n+1}(P,-1)-U_{n-1}(P,-1).  \label{0.6}
\end{equation}

If $P$ is even, then $V_{n}$ is even, $2|U_{n}$ if and only if $2|n,$ and%
\begin{equation}
U_{n}(P,-1)-U_{n-1}(P,-1)\text{ is odd.}  \label{0.10}
\end{equation}

Let $d$ be a positive integer which is not a perfect square and $N$ be any
nonzero fixed integer. Then the equation $x^{2}-dy^{2}=N$ is known as Pell
equation. For $N=1$, the equation $x^{2}-dy^{2}=1$ is known as classical
Pell equation. We use the notations $(x,y)$, and $x+y\sqrt{d}$
interchangeably to denote solutions of the equation $x^{2}-dy^{2}=N.$ Also,
if $x$ and $y$ are both positive, we say that $x+y\sqrt{d}$ is positive
solution to the equation $x^{2}-dy^{2}=N.$ The least positive integer
solution $x_{1}+y_{1}\sqrt{d}$\ to the equation $x^{2}-dy^{2}=N$ is called
the fundamental solution to this equation.

Now, consider the Pell equation 
\begin{equation}
x^{2}-dy^{2}=1.  \label{0.2}
\end{equation}%
If $x_{1}+y_{1}\sqrt{d}~$is the fundamental solution of the equation (\ref%
{0.2}), then all positive integer solutions of this equation are given by 
\[
x_{n}+y_{n}\sqrt{d}=\left( x_{1}+y_{1}\sqrt{d}\right) ^{n} 
\]%
with $n\geq 1.$

\begin{lemma}
\label{L2}Let $x_{1}+y_{1}\sqrt{d}~$be the fundamental solution of the
equation $x^{2}-dy^{2}=1.$ Then all positive integer solutions of the
equation $x^{2}-dy^{2}=1$ are given by 
\[
x_{n}=\frac{V_{n}(2x_{1},-1)}{2}\text{ and }y_{n}=y_{1}U_{n}(2x_{1}-1) 
\]%
with $n\geq 1.$
\end{lemma}

By Theorem 110 given in \cite{NAGEL}, we can give the following lemma.

\begin{lemma}
\label{L12}Let $k_{1}+t_{1}\sqrt{d}$ be fundamental solution of the equation 
$u^{2}-dv^{2}=2$ and let $x_{1}+y_{1}\sqrt{d}$ be the fundamental solution
of the equation $x^{2}-dy^{2}=1.$ Then all positive integer solutions of the
equation $u^{2}-dv^{2}=2$ are given by 
\begin{equation}
X_{n}+Y_{n}\sqrt{d}=\left( k_{1}+t_{1}\sqrt{d}\right) \left( x_{1}+y_{1}%
\sqrt{d}\right) ^{n}  \label{w}
\end{equation}%
with $n\geq 1.$
\end{lemma}

\begin{lemma}
\label{L1}\emph{(\cite{Luca})}If $\ P$ is even, then 
\[
v_{2}(V_{n}(P,-1))=\left\{ 
\begin{array}{c}
v_{2}(P)\text{ \ \ \ }n\equiv 1(\func{mod}2), \\ 
1\text{ \ \ \ \ \ }n\equiv 0(\func{mod}2).%
\end{array}%
\right. 
\]
\end{lemma}

\begin{lemma}
\label{L4}\emph{(\cite{zfr2})}Let $a$ be a positive integer which is not a
perfect square and let $b$ be a positive integer. Let $u_{1}\sqrt{a}+v_{1}%
\sqrt{b}$ be the minimal solution of the equation $ax^{2}-by^{2}=1$ and $%
P=4au_{1}^{2}-2.$ Then all positive integer solutions of the equation $%
ax^{2}-by^{2}=1$ are given by $%
(x,y)=(u_{1}(U_{m+1}-U_{m}),v_{1}(U_{m+1}+U_{m}))$ with $n\geq 0,$ where $%
U_{m}=U_{m}(P,-1).$
\end{lemma}

\begin{lemma}
\label{L10}\emph{(\cite{Zfr})}Let $n\in 
\mathbb{N}
\cup \{0\},$ $m,r\in $ $%
\mathbb{Z}
$ and $m$ be a nonzero integer. Then 
\begin{equation}
U_{2mn+r}(P,-1)\equiv U_{r}(P,-1)(\func{mod}U_{m}(P,-1))  \label{0.8}
\end{equation}%
and 
\[
V_{2mn+r}(P,-1)\equiv V_{r}(P,-1)(\func{mod}U_{m}(P,-1)). 
\]%
From the identity (\ref{0.8}), we can give the followings:%
\begin{equation}
3|U_{n}(P,-1)\Longleftrightarrow \left\{ 
\begin{array}{c}
3|P\text{ and }2|n, \\ 
\text{ \ \ \ \ \ }P\equiv 1(\func{mod}3)\text{ and }3|n,%
\end{array}%
\right.  \label{ax}
\end{equation}%
and
\end{lemma}

\begin{equation}
5|U_{n}(P,-1)\Longleftrightarrow \left\{ 
\begin{array}{c}
P\equiv 0(\func{mod}5)\text{ and }2|n, \\ 
\text{ \ \ \ \ \ }P^{2}\equiv 1(\func{mod}5)\text{ and }3|n, \\ 
P^{2}\equiv -1(\func{mod}5)\text{ and }5|n.%
\end{array}%
\right.  \label{bx}
\end{equation}

\begin{lemma}
\label{L0}Let $P\equiv 1(\func{mod}a)$ and $n\geq 1$. Then $%
a|U_{n}(P,-1)-U_{n-1}(P,-1)$ if and only if $n\equiv 2$ or $5(\func{mod}6).$
\end{lemma}

Throughout this paper, we denote $v_{2}(m)$ by the exponent of $2$ in the
factorization of $m,$ where $m$ is a nonzero positive integer.

\section{Main Theorems}

From now on, unless otherwise indicated, we assume that $m$ and $n$ are
positive integers such that $0<m<n$ and $a,b$ are different fixed positive
integers.

\begin{theorem}
\label{t9}Let $a=2^{t}r$ with $r$ odd and $t\geq 1,$ and let $b$ be odd.
Then the equation $(a^{n}-2^{m})(b^{n}-2^{m})=x^{2}$ has no solutions in the
following cases:
\end{theorem}

$i)$ $2\nmid m.$

$ii)~2|n$ and $2|m.$

$iii)~2\nmid n,~2|m,$~$b\equiv 1(\func{mod}4)$ with $~n-m>1$ or $~n-m=1$ and 
$t>1.$

$iv)~2\nmid n,~2|m,~nt-m=1,$ and $b\equiv 3(\func{mod}4).$

\proof%
Let $a=2^{t}r$ with $r$ odd and $t\geq 1.$ Then we have the equation 
\begin{equation}
(2^{nt-m}r^{n}-1)(b^{n}-2^{m})=\left( \dfrac{x}{2^{m/2}}\right) ^{2}.
\label{12}
\end{equation}

$i)$ Let $2\nmid m.$ It can be seen that the number $v_{2}\left(
2^{m}(2^{nt-m}r^{n}-1)(b^{n}-2^{m})\right) $ is odd. This contradicts the
fact that the number $v_{2}\left( x^{2}\right) $ is even.

$ii)$ Let $2|n$ and $2|m.$ Then from the equation (\ref{12}), since $%
r,~2^{nt-m}r^{n}-1,~$and $b^{n}-2^{m}$ are odd, we get 
\[
\left( \frac{x}{2^{m/2}}\right) ^{2}\equiv (-1)(1)\equiv 3(\func{mod}4), 
\]%
which is impossible.

$iii)~$Let $2\nmid n,~2|m,~$and $b\equiv 1(\func{mod}4).$ If $n-m>1$ or $%
n-m=1$ and $t>1,$ then, in both cases, we get $nt-m>1.$ Thus the equation (%
\ref{12}) yields to 
\[
\left( \dfrac{x}{2^{m/2}}\right) ^{2}\equiv (-1)b\equiv 1(\func{mod}4). 
\]%
This shows that $b\equiv 3(\func{mod}4),$ a contradiction.

$iv)~$Let $2\nmid n,~2|m,~nt-m=1,$ and $b\equiv 3(\func{mod}4).$ If $t>1,$
it yields that $nt>n>m,~$which implies that $nt-m>1.$ Therefore $t=1.$ Then $%
n=m+1=2k+1.$ Thus, by (\ref{12}), we have 
\begin{equation}
2r(r^{k})^{2}-du^{2}=1\text{ and }b(b^{k})^{2}-dv^{2}=2^{2k}  \label{13}
\end{equation}%
for some integers $u$ and $v$ with $\gcd (u,v)=1.$ From the first equation
of (\ref{13}), we get $d\equiv 1(\func{mod}4),$ and thus from the second
equation of (\ref{13}), we get $b\equiv 1(\func{mod}4).$ This is a
contradiction.%
\endproof%

The following lemma can be found in \cite{bhr}.

\begin{lemma}
\label{L7}Let $k\geq 1.$ Then all nonnegative integer solutions of the
equation $u^{2}-5v^{2}=-4^{k}$ are given by $%
(u,v)=(2^{k-1}L_{2m+1},2^{k-1}F_{2m+1})$ with $m\geq 0.$
\end{lemma}

In \cite{szly} , it is shown that the equation $(2^{n}-1)(5^{n}-1)=x^{2}$
has only the solution $(n,x)=(1,2).$

\begin{theorem}
The equation $(2^{n}-2^{m})(5^{n}-2^{m})=x^{2}$ has only the solution $%
(n,m,x)=(3,2,22).$
\end{theorem}

\proof%
Since $1\leq m<n,$ we have $2\nmid n,~2|m,$ and $n-m=1$ by part $iii)$ of
Theorem \ref{t9}$.$ Let $m=2k$ and thus $n=2k+1$ with $k>0.$ In this case,
we get the equation $(5^{2k+1}-2^{2k})=\left( x/2^{k}\right) ^{2},$ i.e., $%
\left( x/2^{k}\right) ^{2}-5(5^{k})^{2}=-4^{k}.$ By Lemma \ref{L7}, it
follows that $x=2^{2k-1}L_{2r+1}$ and $5^{k}=2^{k-1}F_{2r+1}$ with $r\geq 0.$
This implies that $k=1,~r=2$ since $k>0,$ and therefore $n=3,~m=2,$ and $%
x=22.$%
\endproof%

\begin{corollary}
If $(a,b)=(4^{k},5^{t})$ with $k,t\geq 1,$ then the equation $%
(a^{n}-2^{m})(b^{n}-2^{m})=x^{2}$ has no solutions.
\end{corollary}

\begin{theorem}
\label{t11}Let $a$ and $b$ be even integers with $v_{2}(a)\neq v_{2}(b).$
Then the equation $(a^{n}-2^{m})(b^{n}-2^{m})=x^{2}$ has no solutions if $%
2|n.$
\end{theorem}

\proof%
Let $a=2^{t}r$ and $b=2^{l}s$ with $t\neq l$ and $r,s$ odd. Then it can be
seen that 
\begin{equation}
(2^{nt-m}r^{n}-1)(2^{nl-m}s^{n}-1)=\left( x/2^{m}\right) ^{2}.  \label{1.4}
\end{equation}%
Now, we divide the proof into two cases.

Let $m$ be even. In this case, we have the equations 
\[
\left( 2^{\tfrac{nt-m}{2}}r^{\tfrac{n}{2}}\right) ^{2}-du^{2}=1\text{ and }%
\left( 2^{\tfrac{nl-m}{2}}s^{\tfrac{n}{2}}\right) ^{2}-dv^{2}=1 
\]%
for some integers $u$ and $v$ with $\gcd (u,v)=1,$ where $d=\gcd \left(
2^{nt-m}r^{n}-1,2^{nl-m}s^{n}-1\right) .$ Thus by Lemma \ref{L2}, it follows
that 
\[
2^{\tfrac{nt-m}{2}}r^{\tfrac{n}{2}}=\frac{V_{c}(2x_{1},-1)}{2}\text{ and }2^{%
\tfrac{nl-m}{2}}s^{\tfrac{n}{2}}=\frac{V_{k}(2x_{1},-1)}{2}, 
\]%
where $x_{1}+y_{1}\sqrt{d}~$is the fundamental solution of the equation $%
x^{2}-dy^{2}=1.$ Since both of $V_{c}/2$ and $V_{k}/2$ are even, $c$ and $k$
must be odd $\ $and 
\[
v_{2}\left( 2x_{1}\right) =v_{2}\left( 2^{\tfrac{nt-m}{2}+1}r^{\tfrac{n}{2}%
}\right) =v_{2}\left( 2^{\tfrac{nl-m}{2}+1}s^{\tfrac{n}{2}}\right) 
\]%
by Lemma \ref{L1}. From here, we get $\dfrac{nt-m}{2}=\dfrac{nl-m}{2},$
i.e., $t=l,$ which contradicts the fact that $t\neq l.$

Let $m$ be odd. Then, from the equation (\ref{1.4}), we obtain 
\[
2\left( 2^{\tfrac{nt-m-1}{2}}r^{\tfrac{n}{2}}\right) ^{2}-du^{2}=1\text{ and 
}2\left( 2^{\tfrac{nl-m-1}{2}}s^{\tfrac{n}{2}}\right) ^{2}-dv^{2}=1 
\]%
for some integers $u$ and $v$ with $\gcd (u,v)=1,$ where $d=\left(
2^{nt-m}r^{n}-1,2^{nl-m}s^{n}-1\right) .$ By Lemma \ref{L4}, it follows that%
\[
2^{\tfrac{nt-m-1}{2}}r^{\tfrac{n}{2}}=u_{1}(U_{c+1}-U_{c})\text{ and }2^{%
\tfrac{nl-m-1}{2}}s^{\tfrac{n}{2}}=u_{1}(U_{k+1}-U_{k}), 
\]%
where $u_{1}\sqrt{2}+v_{1}\sqrt{d}$ is the minimal solution of the equation $%
2x^{2}-dy^{2}=1$ and $U_{g}=U_{g}(P,-1)$ with $P=8u_{1}^{2}-2.$ Since $%
U_{c+1}-U_{c}$ and $U_{k+1}-U_{k}$ is odd by (\ref{0.10}), it is seen that 
\[
v_{2}\left( u_{1}\right) =v_{2}\left( 2^{\tfrac{nt-m-1}{2}}\right)
=v_{2}\left( 2^{\tfrac{nl-m-1}{2}}\right) . 
\]
This shows that $t=l,$ which contradicts the fact that $t\neq l.$%
\endproof%

\begin{theorem}
\label{t12}Let $a=2^{t}r$ and $b=2^{l}r$ with $t\neq l$ and $r$ odd. Let $t$
and $l$ have the same parity. Then the equation $%
(a^{n}-2^{m})(b^{n}-2^{m})=x^{2}$ has no solutions if $2\nmid n.$
\end{theorem}

\proof%
Since $a=2^{t}r$ and $b=2^{l}r$ with $t\neq l,$ we get 
\[
(2^{nt-m}r^{n}-1)(2^{nl-m}r^{n}-1)=\left( x/2^{m}\right) ^{2}. 
\]%
Let $n=2k+1.$ We shall discuss separately the proof according to whether $m$
is even or odd.

Let $m$ be even. According to whether $t$ and $l$ are even or odd, we obtain
the following equations, respectively, 
\begin{equation}
r\left( 2^{\tfrac{nt-m}{2}}r^{k}\right) ^{2}-du^{2}=1\text{ and }r\left( 2^{%
\tfrac{nl-m}{2}}r^{k}\right) ^{2}-dv^{2}=1  \label{x}
\end{equation}%
or 
\begin{equation}
2r\left( 2^{\tfrac{nt-m-1}{2}}r^{k}\right) ^{2}-du^{2}=1\text{ and }2r\left(
2^{\tfrac{nl-m-1}{2}}r^{k}\right) ^{2}-dv^{2}=1  \label{y}
\end{equation}%
for some integers $u$ and $v$ with $\gcd (u,v)=1,$ where $d=\gcd \left(
2^{nt-m}r^{n}-1,2^{nl-m}r^{n}-1\right) .$ Lemma \ref{L4} and (\ref{x})\
implies that if $r$ is not a perfect square, then%
\begin{equation}
2^{\tfrac{nt-m}{2}}r^{k}=u_{1}(U_{c+1}-U_{c})\text{ and }2^{\tfrac{nl-m}{2}%
}r^{k}=u_{1}(U_{k+1}-U_{k}),  \label{a}
\end{equation}%
and if $r$ is a perfect square, then 
\begin{equation}
2^{\tfrac{nt-m}{2}}z^{2k+1}=\dfrac{V_{s}(2x_{1},-1)}{2}\text{ and }2^{\tfrac{%
nl-m}{2}}z^{2k+1}=\dfrac{V_{g}(2x_{1},-1)}{2}.  \label{b}
\end{equation}%
where $u_{1}\sqrt{r}+v_{1}\sqrt{d}$ and $x_{1}+y_{1}\sqrt{d}$ are,
respectively, the minimal solution of the equations $rx^{2}-dy^{2}=1$ and $%
~x^{2}-dy^{2}=1.$ Similarly, (\ref{y}) implies that if $r$ is not a perfect
square, then 
\begin{equation}
2^{\tfrac{nt-m-1}{2}}r^{k}=u_{2}(U_{c+1}-U_{c})\text{ and }2^{\tfrac{nl-m-1}{%
2}}r^{k}=u_{2}(U_{k+1}-U_{k}),  \label{c}
\end{equation}%
and if $r=z^{2},$ then 
\begin{equation}
2^{\tfrac{nt-m-1}{2}}z^{2k+1}=u_{3}(U_{c+1}-U_{c})\text{ and }2^{\tfrac{%
nl-m-1}{2}}z^{2k+1}=u_{3}(U_{k+1}-U_{k}).  \label{d}
\end{equation}%
where $u_{2}\sqrt{2r}+v_{2}\sqrt{d}$ and $u_{3}\sqrt{2}+v_{3}\sqrt{d}$ are
the minimal solution of the equations $2rx^{2}-dy^{2}=1,$ and $%
2x^{2}-dy^{2}=1,$ respectively. Since $U_{c+1}-U_{c}$ and $U_{k+1}-U_{k}$ is
odd by (\ref{0.10}), it is seen from (\ref{a}), (\ref{c}), and (\ref{d})
that 
\[
v_{2}\left( u_{1}\right) =v_{2}\left( 2^{\tfrac{nt-m}{2}}\right)
=v_{2}\left( 2^{\tfrac{nl-m}{2}}\right) , 
\]%
or 
\[
v_{2}\left( u_{2}\right) =v_{2}\left( 2^{\tfrac{nt-m-1}{2}}\right)
=v_{2}\left( 2^{\tfrac{nl-m-1}{2}}\right) , 
\]%
or%
\[
v_{2}\left( u_{3}\right) =v_{2}\left( 2^{\tfrac{nt-m-1}{2}}\right)
=v_{2}\left( 2^{\tfrac{nl-m-1}{2}}\right) . 
\]%
In all the above cases, we get $t=l,$ which contradicts the fact that $t\neq
l.$ Since $V_{s}(2x_{1},-1)/2$ and $\ V_{g}(2x_{1},-1)/2$ are even,$\ s$ and 
$g$ are odd\ and 
\[
v_{2}\left( 2x_{1}\right) =v_{2}\left( 2^{\tfrac{nt-m}{2}}z^{2k+1}\right)
=v_{2}\left( 2^{\tfrac{nl-m}{2}}z^{2k+1}\right) 
\]%
by Lemma \ref{L1}. From here, we get $\dfrac{nt-m}{2}=\dfrac{nl-m}{2},$
i.e., $t=l,$ which contradicts the fact that $t\neq l.$

A similar argument shows that the case $m$ odd is also impossible.%
\endproof%

From Theorems \ref{t11} and \ref{t12}, we can conclude the following result.

\begin{corollary}
\label{c14}Let $a$ be even, $b=2^{k}a$ with $k\geq 1$ and let $v_{2}(a)$ and 
$v_{2}(b)$ have the same parity. Then the equation $%
(a^{n}-2^{m})(b^{n}-2^{m})=x^{2}$ has no solutions for all $m.$
\end{corollary}

\begin{theorem}
\label{t13}Let $a\equiv 2(\func{mod}3)$ and $3|b.$ Then the equation $%
(a^{n}-2^{m})(b^{n}-2^{m})=x^{2}$ has no solutions in the following cases:
\end{theorem}

$i)$ If $m$ and $n$ have opposite parity.

$ii)$ If $2|n,2|m,$ and $a$ and $b$ are both even.

\proof%
$i)$ Assume that $m$ and $n$ have opposite parity. Then it can be seen that $%
x^{2}=(a^{n}-2^{m})(b^{n}-2^{m})\equiv 2(\func{mod}3),$ which is impossible.

$ii)$ Assume that $2|n,2|m,$ and $a$ and $b$ are both even. Let $a=2^{t}r$
and $b=2^{l}s$ with $r,s$ odd. Then it can be seen that 
\begin{equation}
(2^{nt-m}r^{n}-1)(2^{nl-m}s^{n}-1)=\left( x/2^{m}\right) ^{2}.  \label{1.6}
\end{equation}%
Thus we get 
\[
\left( 2^{\tfrac{nt-m}{2}}r^{\tfrac{n}{2}}\right) ^{2}-du^{2}=1\text{ and }%
\left( 2^{\tfrac{nl-m}{2}}s^{\tfrac{n}{2}}\right) ^{2}-dv^{2}=1 
\]%
for some integers $u$ and $v$ with $(u,v)=1,$ where $d=\gcd \left(
2^{nt-m}r^{n}-1,2^{nl-m}s^{n}-1\right) .$ Thus by Lemma \ref{L2}, it follows
that 
\[
2^{\tfrac{nt-m}{2}}r^{\tfrac{n}{2}}=\frac{V_{k_{1}}(2x_{1},-1)}{2}%
,~u=y_{1}U_{k_{1}}(2x_{1},-1)\text{ } 
\]%
and 
\[
2^{\tfrac{nl-m}{2}}s^{\tfrac{n}{2}}=\frac{V_{k_{2}}(2x_{1},-1)}{2}%
,~v=y_{1}U_{k_{2}}(2x_{1},-1), 
\]%
where $x_{1}+y_{1}\sqrt{d}~$is the fundamental solution of the equation $%
x^{2}-dy^{2}=1.$ Without loss of generality, assume that $t\geq l.$ Since $%
\gcd (u,v)=1,$ by (\ref{0.4}) and (\ref{0.5}), we get $y_{1}=1,~\gcd
(k_{1},k_{2})=1,$ and $x_{1}=2^{\tfrac{nl-m}{2}}g^{\tfrac{n}{2}}$, where $%
g=\gcd (r,s)$. Thus $d=2^{nl-m}g^{n}-1.$ This shows that $d\equiv 0(\func{mod%
}3).$ On the other hand, it can be seen that $dv^{2}\equiv d\equiv -1(\func{%
mod}3),$ a contradiction.%
\endproof%

From Theorems \ref{t9} and \ref{t13}, we can deduce the following corollary.

\begin{corollary}
If $a\equiv 2(\func{mod}6)$ and $b\equiv 3(\func{mod}6),$ then the equation $%
(a^{n}-2^{m})(b^{n}-2^{m})=x^{2}$ has no solutions.
\end{corollary}

Thus, taking $m=1$ in the equation (\ref{3.1}), we can give the following
four corollaries from the above theorems.

\begin{corollary}
\label{t1}Let $a$ and $b$ have opposite parity. Then the equation $%
(a^{n}-2)(b^{n}-2)=x^{2}$ has no solutions for $n>1.$
\end{corollary}

\begin{corollary}
\label{t4}Let $3|a$ and $3\nmid b.$ Then the equation $%
(a^{n}-2)(b^{n}-2)=x^{2}$ has no solutions if $2|n.$
\end{corollary}

\begin{corollary}
\label{t2}Let $a$ and $b$ be even integers with $v_{2}(a)\neq v_{2}(b).$
Then the equation $(a^{n}-2)(b^{n}-2)=x^{2}$ has no solutions if $2|n.$
\end{corollary}

\begin{corollary}
\label{t3}Let $a=2^{t}r$ and $b=2^{l}r$ with $t\neq l$ and let $t$ and $l$
have the same parity. Then the equation $(a^{n}-2)(b^{n}-2)=x^{2}$ has no
solutions if $2\nmid n.$
\end{corollary}

Immediately, from last two corollaries, we can conclude the following result.

\begin{corollary}
\label{c1}Let $a=2^{t}r$ and $b=2^{l}r$ with $t\neq l$ and let $t$ and $l$
have the same parity. Then the equation $(a^{n}-2)(b^{n}-2)=x^{2}$ has no
solutions if $n>1.$
\end{corollary}

If $a+b\sqrt{d}$ is a solution of the equation $x^{2}-dy^{2}=2,$ then $(a+b%
\sqrt{d})^{2}/2=$ $\left( a^{2}+db^{2}\right) /2+ab\sqrt{d}$ is a solution
of the equation $x^{2}-dy^{2}=1.$

The proof of the following lemma can be found in \cite{rfk}.

\begin{lemma}
\label{L5}Let $d>2.$ If $k_{1}+t_{1}\sqrt{d}$ is the fundamental solution of
the equation $u^{2}-dv^{2}=2,$ then $\left( k_{1}^{2}+dt_{1}^{2}\right)
/2+k_{1}t_{1}\sqrt{d}$ is the fundamental solution of the equation $%
x^{2}-dy^{2}=1.$
\end{lemma}

\begin{theorem}
\label{t6}Let $d>2.$ Let $x_{1}+y_{1}\sqrt{d}$ and $k_{1}+t_{1}\sqrt{d}$ $~$%
be the fundamental solutions of the equations$\ x^{2}-dy^{2}=1$ and $%
u^{2}-dv^{2}=2,$ respectively. Then $(x_{1},y_{1})=(k_{1}^{2}-1,k_{1}t_{1})$
and all solutions of the equation $u^{2}-dv^{2}=2$ are given by 
\[
(X_{n},Y_{n})=\left( k_{1}\left( U_{n+1}-U_{n}\right) ,t_{1}\left(
U_{n+1}+U_{n}\right) \right) 
\]%
with $n\geq 1,$ where $U_{n}=U_{n}(2x_{1},-1)$.
\end{theorem}

\proof%
Since $k_{1}+t_{1}\sqrt{d}$ $~$be the fundamental solutions of the equation $%
u^{2}-dv^{2}=2,$ $x_{1}+y_{1}\sqrt{d}=\left( k_{1}^{2}+dt_{1}^{2}\right)
/2+k_{1}t_{1}\sqrt{d}$ \ is the fundamental solution of the equation $%
x^{2}-dy^{2}=1$ by Lemma \ref{L5}. It is immediately seen that $%
(x_{1},y_{1})=(k_{1}^{2}-1,k_{1}t_{1}).$ By (\ref{w}), all positive solution
of the equation $u^{2}-dv^{2}=2$ are given by $X_{n}+Y_{n}\sqrt{d}=\left(
k_{1}+t_{1}\sqrt{d}\right) \left( x_{1}+y_{1}\sqrt{d}\right) ^{n}$.
Therefore $X_{n}+Y_{n}\sqrt{d}=\left( k_{1}+t_{1}\sqrt{d}\right) \left(
x_{n}+y_{n}\sqrt{d}\right) .$ It follows that%
\begin{eqnarray*}
X_{n}+Y_{n}\sqrt{d} &=&\left( \frac{k_{1}^{2}+y_{1}\sqrt{d}}{k_{1}}\right)
\left( x_{1}+y_{1}\sqrt{d}\right) ^{n} \\
&=&\left( \frac{x_{1}+1+y_{1}\sqrt{d}}{k_{1}}\right) \left( x_{1}+y_{1}\sqrt{%
d}\right) ^{n} \\
&=&\frac{1}{k_{1}}\left[ \left( x_{1}+y_{1}\sqrt{d}\right) ^{n+1}+\left(
x_{1}+y_{1}\sqrt{d}\right) ^{n}\right] .
\end{eqnarray*}%
Thus 
\[
X_{n}=\dfrac{x_{n+1}+x_{n}}{k_{1}}\text{ and }Y_{n}=\dfrac{y_{n+1}+y_{n}}{%
k_{1}} 
\]%
Since $(x_{n},y_{n})$ is a solution of the equation$~x^{2}-dy^{2}=1,$ $%
(x_{n},y_{n})=\left( V_{n}(2x_{1},-1)/2,y_{1}U_{n}(2x_{1},-1)\right) $ by
Lemma \ref{L2}. Using the identity (\ref{0.6}), we get 
\begin{eqnarray*}
X_{n} &=&\frac{1}{2k_{1}}\left( V_{n+1}+V_{n}\right) =\frac{1}{2k_{1}}\left(
U_{n+2}-U_{n}+U_{n+1}-U_{n-1}\right) \\
&=&\frac{1}{2k_{1}}\left(
2x_{1}U_{n+1}-U_{n}-U_{n}+U_{n+1}+U_{n+1}-PU_{n}\right) \\
&=&\frac{1}{2k_{1}}\left[ \left( 2x_{1}+2\right) U_{n+1}-\left(
2x_{1}+2\right) U_{n}\right] \\
&=&\frac{1}{2k_{1}}\left[ 2k_{1}^{2}U_{n+1}-2k_{1}^{2}U_{n}\right] \\
&=&k_{1}\left( U_{n+1}-U_{n}\right)
\end{eqnarray*}%
and 
\[
Y_{n}=\frac{1}{k_{1}}\left( y_{1}U_{n+1}+y_{1}U_{n}\right) =\frac{y_{1}}{%
k_{1}}\left( U_{n+1}+U_{n}\right) =t_{1}\left( U_{n+1}+U_{n}\right) =.%
\endproof%
\]

Now, we can give the following theorem.

\begin{theorem}
\label{t7}Let $\gcd (a,b)=1.$ Then the equation $(a^{n}-2)(b^{n}-2)=x^{2}$
has no solutions if $2|n.$
\end{theorem}

\proof%
Assume that $\gcd (a,b)=1$ and $2|n.$ Let $n=2m.$ Then 
\begin{equation}
\left( a^{m}\right) ^{2}-du^{2}=2\text{ and }\left( b^{m}\right)
^{2}-dv^{2}=2,  \label{1.3}
\end{equation}%
for some integers $u$ and $v$ with $\gcd (u,v)=1,$ where $%
d=(a^{n}-2,b^{n}-2).$ It is obvious that $d>2.$ Assume that $(x_{1},y_{1})$
and $(k_{1},t_{1})$ are the fundamental solutions of the equations $%
x^{2}-dy^{2}=1$ and $u^{2}-dv^{2}=2,$ respectively. Then by Theorem \ref{t6}%
, we get 
\[
a^{m}=k_{1}(U_{r}-U_{r-1}),\text{ }u=t_{1}(U_{r}+U_{r-1}) 
\]%
and 
\[
b^{m}=k_{1}(U_{s}-U_{s-1}),\text{ }v=t_{1}(U_{s}+U_{s-1}). 
\]
Since $\gcd (a,b)=1$, $\ $we get $\gcd (a^{m},b^{m})=1.$ From the above
equations it follows that $(k_{1},t_{1})=(1,1)$ since $\gcd
(a^{m},b^{m})=\gcd (u,v)=1.$ This implies that $d=-1,$ which is impossible.%
\endproof%

Using Mathematica, we verified for all $2\leq a<b\leq 100$ and $n$ in the
range $2\leq n\leq 1000$ that the equation $(a^{n}-2)(b^{n}-2)=x^{2}$ has
only solutions $%
(a,b,n,x)=(2,10,2,14),(2,10,6,7874),(2,58,2,82),(3,45,2,119),(4,100,3,7874),(10,58,2,574). 
$ Now, we solve the equation $(a^{n}-2)(b^{n}-2)=x^{2}$ for $%
(a,b)=(2,10),(3,45),$ $(4,100),(10,58).$

The proofs of the following two lemmas can be done by induction on $m.$

\begin{lemma}
\label{L9}Let $m\geq 4.$ Then $5^{m}>2^{2m+1}-3.$
\end{lemma}

\begin{lemma}
\label{L11}Let $m\geq 2.$ Then $2\cdot 3^{4m-3}>5^{m}+1.$
\end{lemma}

From Lemma \ref{L9}, we can give the following corollary.

\begin{corollary}
\label{C1}Let $m$ and $z$ be positive integers. If $(z+1)(2z-1)^{2}=10^{2m}$%
, then $m=1,z=3$ or $m=3,z=63.$
\end{corollary}

\begin{theorem}
\label{t15}The equation $(2^{n}-2)(10^{n}-2)=x^{2}$ has only the solutions $%
(n,x)=(2,14),(6,6874).$
\end{theorem}

\proof%
It is obvious that $(n,x)=(2,14)$ is a solution. Let $n>2.$ Firstly, assume
that $n$ is even, say $n=2m.$ Then 
\[
\left( 2^{m}\right) ^{2}-2du^{2}=2\text{ and }\left( 10^{m}\right)
^{2}-2dv^{2}=2\text{ } 
\]%
for some integers $u$ and $v$ with $\gcd (u,v)=1,$ where $%
2d=(2^{n}-2,10^{n}-2).$ Since $m>1$, it can be seen that $2d>2.$\ Hence, by
Theorem \ref{t6}, it follows that 
\[
2^{m}=k_{1}(U_{r}-U_{r-1}),\text{ }u=t_{1}(U_{r}+U_{r-1}) 
\]%
and 
\[
10^{m}=k_{1}(U_{s}-U_{s-1}),\text{ }v=t_{1}(U_{s}+U_{s-1}), 
\]%
where $U_{t}=U_{t}(2x_{1},-1)$ and $(x_{1},y_{1})$ and $(k_{1},t_{1})$ are
the fundamental solutions of the equations $x^{2}-2dy^{2}=1$ and $%
u^{2}-2dv^{2}=2,$ respectively. Since $U_{r}+U_{r-1}$ is odd by (\ref{0.10})
and $\gcd (u,v)=1,$ it follows that $k_{1}=2^{m}$, $r=1,$ and $t_{1}=1,$
which implies that $u=1.$ Thus $\left( 2^{m}\right) ^{2}-2d=2,$ i.e., $%
d=2^{n-1}-1.$ By Lemma \ref{L5}, $x_{1}+y_{1}\sqrt{2d}=(2^{m}+\sqrt{2d}%
)^{2}/2=2^{n}-1+2^{m}\sqrt{2d}.$ This shows that $x_{1}=2^{n}-1$ and $%
y_{1}=2^{m}.$ On the other hand, since $(10^{m}+v\sqrt{2d}%
)^{2}=10^{n}-1+10^{m}v\sqrt{2d}$ is a solution of the equation $%
x^{2}-2dy^{2}=1,$ it follows that 
\[
10^{n}-1=V_{k}(2x_{1},-1)/2\text{ and }%
10^{m}v=y_{1}U_{k}(2x_{1},-1)=2^{m}U_{k} 
\]%
for some positive integer $k$ by Lemma \ref{L2}. This implies that $5|U_{k}.$
Now assume that $4|n.$ Then $U_{2}=P=2x_{1}=2^{n+1}-2\equiv 0(\func{mod}5)$
and this implies that $k$ is even since $5|U_{k}.$ Taking $k=2c,$ we get $%
2\cdot 10^{n}-2=V_{2c}=V_{c}^{2}-2$ by (\ref{0.11}), i.e., $2\cdot
10^{n}=V_{c}^{2}.$ This is impossible. Hence $2||n.$ Then $P\equiv 1(\func{%
mod}5).$ Since $U_{k}=U_{6q+r}\equiv U_{r}(\func{mod}U_{3})$ by (\ref{0.8})
and $5|U_{3},$ it follows that $k=3t.$ This implies that $2\cdot
10^{n}-2=V_{3t}=V_{t}^{3}-3V_{t}$ by (\ref{0.9}). Taking $V_{t}=2z,$ then,
from the last equality, we get $(z+1)(2z-1)^{2}=10^{2m}.$ By Corollary \ref%
{C1}, it follows that $m=3.$ This shows that $n=6.$

Secondly assume that $n$ is odd. If $n=4k+3,$ this gives $x^{2}\equiv 3(%
\func{mod}5),$ a contradiction. If $n=4k+1,$ then we have the equation $%
(2^{4k}-1)(2^{4k}5^{4k+1}-1)=\left( \dfrac{x}{2}\right) ^{2}.$ This implies
that 
\[
\left( 2^{k}\right) ^{4}-du^{2}=1\text{ and }5\left( 10^{k}\right)
^{4}-dv^{2}=1 
\]%
for some integers $u$ and $v$ with $\gcd (u,v)=1,$ where $d=\gcd
(2^{4k}-1,2^{4k}5^{4k+1}-1).$ From the equation $\left( 2^{k}\right)
^{4}-du^{2}=1$, we can write that $2^{2k}-1=ra^{2},2^{2k}+1=sb^{2}$ for some
integers $a,b,r,s,$ where $rs=d.$ The equation $5\left( 10^{k}\right)
^{4}-dv^{2}=1$ implies that $5\left( 10^{2k}\right) ^{2}\equiv 1(\func{mod}%
r) $ and $5\left( 10^{2k}\right) ^{2}\equiv 1(\func{mod}s).$ Thus $\left( 
\dfrac{r}{5}\right) =\left( \dfrac{s}{5}\right) =1.$ If $k$ is even, then $%
sb^{2}\equiv 2(\func{mod}5),$ which implies that $\left( \dfrac{s}{5}\right)
=-1.$ This is a contradiction. If $k$ is odd, then $ra^{2}\equiv 3(\func{mod}%
5).$ This shows that $\left( \dfrac{r}{5}\right) =\left( \dfrac{3}{5}\right)
-1,$ which is a contradiction. This completes the proof.%
\endproof%

From Theorem \ref{t15}, we immediately deduce the following corollary.

\begin{corollary}
\label{C11}The equation $(4^{n}-2)(100^{n}-2)=x^{2}$ has only the solution $%
(n,x)=(3,7874).$
\end{corollary}

\begin{theorem}
\label{t16}The equation $(10^{n}-2)(58^{n}-2)=x^{2}$ has only the solution $%
(n,x)=(2,574).$
\end{theorem}

\proof%
If $4|n$ or $n\equiv 1(\func{mod}4),$ then it can be seen that $x^{2}\equiv
2 $ or$\ 3(\func{mod}5),$ which is impossible. Assume that $n\equiv 3(\func{%
mod}4).$ Then $n$ is form of \ $12q+3$ or $12q+7,$ or $12q+11.$ If $n=12q+3,$
then we get $x^{2}\equiv 3(\func{mod}7),$ a contradiction. If $n=12q+11,$
then it is seen that $x^{2}\equiv 6(\func{mod}7),$ which is impossible. Let $%
n=12q+7.$ Then we get $n\equiv 7,19,31,43,$ or $55(\func{mod}60).$ If $%
n\equiv 31(\func{mod}60),$ then we get $x^{2}\equiv 8(\func{mod}11),$ a
contradiction. Let $n\equiv 7,19,43,$ or $55(\func{mod}60).$ Similarly, when
we investigate the equation $(10^{n}-2)(58^{n}-2)=x^{2}$ according to modulo 
$31,$ we can see that it has no solutions.

Now assume that $n\equiv 2(\func{mod}4).$ Say $n=2m$ with $m$ odd. Then we
get $(2^{2m-1}5^{2m}-1)(2^{2m-1}29^{2m}-1)=x^{2}.$ Thus 
\[
2\left( 2^{m-1}5^{m}\right) ^{2}-du^{2}=1\text{ and }2\left(
2^{m-1}29^{m}\right) ^{2}-dv^{2}=1\text{ } 
\]%
for some integers $u$ and $v$ with $\gcd (u,v)=1,$ where $d=\gcd
(2^{2m-1}5^{2m}-1,2^{2m-1}29^{2m}-1).$ By Lemma \ref{L4}, we obtain 
\[
2^{m-1}5^{m}=u_{1}(U_{r}-U_{r-1}),\text{ }u=v_{1}(U_{r}+U_{r-1}) 
\]%
and 
\[
2^{m-1}29^{m}=u_{1}(U_{s}-U_{s-1}),\text{ }v=v_{1}(U_{s}+U_{s-1}), 
\]%
where $U_{c}=U_{c}(P,-1)$ with $P=8u_{1}^{2}-2,$ and $(u_{1},v_{1})$ is the
fundamental solution of the equation $2x^{2}-dy^{2}=1.$ Since $U_{r}+U_{r-1}$
and $U_{s}+U_{s-1}$ are odd by (\ref{0.10}) and $\gcd (u,v)=1,$ it follows
that $u_{1}=2^{m-1},$ $v_{1}=1,U_{r}-U_{r-1}=5^{m}$ and $%
U_{s}-U_{s-1}=29^{m}.$ Thus $2\left( 2^{m-1}\right) ^{2}-d=1,$ i.e., $%
d=2^{n-1}-1.$ Since $2^{m-1}\sqrt{2}+\sqrt{d}$ is the fundamental solution
of the equation $2x^{2}-dy^{2}=1,$ $\left( 2^{m-1}\sqrt{2}+\sqrt{d}\right)
^{2}=2^{n}-1+2^{m}\sqrt{2d}$ is the fundamental solution of the equation $%
x^{2}-2dy^{2}=1.$ On the other hand, it is seen that $\left( 2^{m-1}5^{m}%
\sqrt{2}+u\sqrt{d}\right) ^{2}=10^{n}-1+10^{m}u\sqrt{2d},$ and this is a
solution of the equation $x^{2}-2dy^{2}=1.$ Hence, by Lemma \ref{L2}, we get 
$10^{n}-1=V_{k}(2x_{1},-1)/2$ and $10^{m}u=y_{1}U_{k}(2x_{1},-1),$ where $%
x_{1}+y_{1}\sqrt{2d}=2^{n}-1+2^{m}\sqrt{d}.$ From this, it is clear that $%
V_{k}(2x_{1},-1)=2\cdot 10^{n}-2$ and $U_{k}(2x_{1},-1)=5^{m}u.$ Therefore $%
5|U_{k}.$ Since $2x_{1}=2^{n+1}-2\equiv 1(\func{mod}5),$ it follows that $%
3|k $ by (\ref{bx}) $.$ Let $k=3t.$ Then we get $2\cdot
10^{n}-2=V_{3t}=V_{t}^{3}-3V_{t}$ by (\ref{0.9}). Taking $V_{t}=2z,$ then
from the last equality we get $(z+1)(2z-1)^{2}=10^{2m}.$ By Corollary \ref%
{C1}, it follows that $m=1$ or $m=3.$ Therefore $m=1$ or $m=3.$ This shows
that $n=2$ or $n=6.$ But $n=6$ is impossible. This completes the proof.%
\endproof%

\begin{theorem}
\label{t17}The equation $(3^{n}-2)(45^{n}-2)=x^{2}$ has only the solution $%
(n,x)=(2,119).$
\end{theorem}

\proof%
If $4|n$ or $n\equiv 1(\func{mod}4),$ then it can be seen that $x^{2}\equiv
2 $ or$\ 3(\func{mod}5),$ which is impossible. Assume that $n\equiv 3(\func{%
mod}4).$ Then $n\equiv 3,7,$ or $11(\func{mod}12).$ In these cases, if using
modulo $13,$ it can be seen that the equation $(3^{n}-2)(45^{n}-2)=x^{2}$ is
impossible.

Now assume that $n\equiv 2(\func{mod}4).$ Say $n=2m$ with $m$ odd$.$ Then we
get 
\[
\left( 3^{m}\right) ^{2}-du^{2}=2\text{ and }\left( 45^{m}\right)
^{2}-dv^{2}=2\text{ } 
\]%
for some integers $u$ and $v$ with $\gcd (u,v)=1,$ where $%
d=(3^{n}-2,45^{n}-2).$ By Theorem \ref{t6}, we obtain 
\begin{equation}
3^{m}=k_{1}(U_{r}-U_{r-1}),\text{ }u=t_{1}(U_{r}+U_{r-1})  \label{ay}
\end{equation}%
and 
\begin{equation}
45^{m}=k_{1}(U_{s}-U_{s-1}),\text{ }v=t_{1}(U_{s}+U_{s-1}),  \label{by}
\end{equation}%
where $U_{c}=U_{c}(P,-1)$ with $P=2x_{1}=2k_{1}^{2}-2,~y_{1}=k_{1}t_{1},$ $%
(x_{1},y_{1})$ is the fundamental solution of the equation $x^{2}-dy^{2}=1,$
and $(k_{1},t_{1})$ is the fundamental solution of the equation $%
x^{2}-dy^{2}=2$. So $k_{1}>1$ and this implies that $3|k_{1}$ by (\ref{ay}).
Then $P=2k_{1}^{2}-2\equiv 1(\func{mod}3).$ Assume that $3|U_{r}-U_{r-1}$.
Then it is seen that $r\equiv 2,5(\func{mod}6)$ by Lemma \ref{L0}. Since $%
5\nmid k_{1}$ by (\ref{ay}), it follows that $P=2k_{1}^{2}-2\equiv 0,1(\func{%
mod}5).$ Assume that $P\equiv 1(\func{mod}5).$ Then $5|U_{r}-U_{r-1}$ by
Lemma \ref{L0} since $r\equiv 2,5(\func{mod}6).$ This is impossible by (\ref%
{ay}). \ Let $P\equiv 0(\func{mod}5).$ Then it can be seen that $5\nmid
U_{s}-U_{s-1}$ by (\ref{bx})$.$ This is impossible by (\ref{by}) since $%
5\nmid k_{1}.$ We conclude that $3\nmid U_{r}-U_{r-1}$ and therefore $%
k_{1}=3^{m}$ and $r=0$ by (\ref{ay}). Thus $x_{1}=k_{1}^{2}-1=3^{n}-1.$
Besides, it is clear that $u=t_{1}=1$ since $\gcd (u,v)=1$ and $r=0.$ This
implies that $y_{1}=k_{1}t_{1}=3^{m}$ and $d=3^{n}-2.$ On the other hand,
since $45^{m}+v\sqrt{d}$ is a solution of the equation $x^{2}-dy^{2}=2,$ $%
\left( 45^{m}+v\sqrt{d}\right) ^{2}/2=45^{n}-1+45^{m}v\sqrt{d}$ is a
solution of the equation $x^{2}-dy^{2}=1$. In this case, $%
45^{n}-1=V_{k}(P,-1)/2$ and $45^{m}v=y_{1}U_{k}(P,-1)=3^{m}U_{k}(P,-1)$ by
Lemma \ref{L2}. From this, it is clear that $V_{k}(P,-1)=2\cdot 45^{2m}-2$
and $U_{k}(P,-1)=3^{m}5^{m}v.$ Therefore $3|U_{k}.$ Since $P\equiv 1(\func{%
mod}3),$ it follows that $3|k$ by (\ref{ax})$.$ Let $k=3t.$ Then we get $%
2\cdot 45^{2m}-2=V_{3t}=V_{t}^{3}-3V_{t}$ by (\ref{0.9}). Taking $V_{t}=2z,$
from the last equality, we get $(z+1)(2z-1)^{2}=w^{2}$ with $w=45^{m}.$ Then 
$z+1=(\frac{w}{2z-1})^{2}.$ Let $a=\frac{w}{2z-1}.$ Then $z+1=a^{2}$ and $%
w=a(2z-1)=a(2a^{2}-3).$ It can be seen that $\gcd (a,2a^{2}-3)=3$ since $%
3|a(2a^{2}-3).$ Then $a=3b$ and $2a^{2}-3=3c$ for some integers $b$ and $c$
with $\gcd (b,c)=1.$ Then $9bc=45^{m}=3^{2m}5^{m}$ and therefore b$%
c=3^{2m-2}5^{m}.$ Since $a=3b$ and $2a^{2}-3=3c,$ we get $c=6b^{2}-1$ and
therefore $b(6b^{2}-1)=3^{2m-2}5^{m}.$ Then $b=3^{2m-2}$ and $6b^{2}-1=5^{m}$%
. Thus $2\cdot 3^{4m-3}=5^{m}+1.$ This is only possible for $m=1$ by Lemma %
\ref{L11}. Consequently, $n=2$ and $x=119.$%
\endproof%

Now we present a conjecture which we are not able to prove.

\begin{conjecture}
Let $k>3$ be odd. Then the equation $(2^{n}-2)((2P_{k})^{n}-2)=x^{2}$ has
only the solution $(n,x)=(2,Q_{k})$.
\end{conjecture}

When $k=3,$ $P_{k}=P_{3}=5$ and the above equation becomes $%
(2^{n}-2)((10^{n}-2)=x^{2},$ which has only the solutions $%
(n,x)=(2,14),(6,7874)$ by Theorem \ref{t15}.

We think the following conjecture is true.

\begin{conjecture}
Let $2<a<b.$ If the equation $(a^{n}-2)(b^{n}-2)=x^{2}$ has a solution $n,x,$
then $n\leq 6.$
\end{conjecture}

\end{document}